\documentclass{amsart}
\usepackage[latin1]{inputenc}
\usepackage{amssymb}
\usepackage[usenames,dvipsnames,dvips]{color}
\usepackage{xcolor}

\newtheorem{thm}{Theorem}
\newtheorem{cor}[thm]{Corollary}
\newtheorem{lem}[thm]{Lemma}
\newtheorem{prop}[thm]{Proposition}
\newtheorem{problem}[thm]{Problem}
\theoremstyle{definition}
\newtheorem{defn}[thm]{Definition}
\theoremstyle{remark}
\newtheorem{rem}[thm]{Remark}
\numberwithin{equation}{section}
\newcommand{\To}{\longrightarrow}

\begin{document}

%\title{Polynomials and multilinear forms in Banach spaces}

\title{Zero subspaces of polynomials on $\ell_1(\Gamma)$}

\author{Antonio Avil\'{e}s and Stevo Todorcevic}


\address{A. Avil\'{e}s, University of Paris 7, UMR 7056, 2, Place Jussieu, Case 7012, 75251 Paris (France)} \email{avileslo@um.es}
\address{S. Todorcevic, C.N.R.S - University of Paris 7, UMR 7056, 2, Place Jussieu, Case 7012, 75251 Paris (France) and Department of Mathematics, University of Toronto, Toronto, M5S 3G3 (Canada)}
\email{stevo@math.jussieu.fr, stevo@math.toronto.edu}

\thanks{The first author was supported by a Marie Curie
Intra-European Felloship MCEIF-CT2006-038768 (E.U.) and research
projects MTM2005-08379 and S\'{e}neca 00690/PI/04 (Spain)}

\begin{abstract}
We provide two examples of complex homogeneous quadratic polynomials
$P$ on Banach spaces of the form $\ell_1(\Gamma)$. The first
polynomial $P$ has both separable and nonseparable maximal zero
subspaces. The second polynomial $P$
has the property that while the index-set $\Gamma$ is not countable, all zero subspaces of $P$ are separable.
\end{abstract}

\maketitle

\section{Introduction}

Let $X$ be a Banach space. A multilinear form on $X$ of degree $n$
is a multilinear function $A:X^n\To Z$ which is moreover
continuous. Continuity is equivalent to the existence of a
constant $K$ such that $\|A(x_1,\ldots,x_n)\|\leq
K\|x_1\|\cdots\|x_n\|$. The minimal such $K$ is the norm of $A$.
The multilinear form $A$ is said to be symmetric if
$A(x_1,\ldots,x_n) = A(x_{\sigma(1)},\ldots,x_{\sigma(n)})$ for
every permutation $\sigma$. A homogenous polynomial of degree $n$
on $X$ is a function $P:X\To Z$ of the form $P(x)=A(x,\ldots,x)$
for $A$ a symmetric multilinear form on $X$. This multilinear form
$A$ is uniquely determined by $P$. In general, a polynomial is a
sum of homogenuous polynomials and constants (constant functions
are considered as homogenous polynomials of degree 0). A
homogeneous polynomial of degree two is called a quadratic
functional.\\

A zero subspace of a polynomial $P$ is a subspace $Y\subset X$ such
that $P(y)=0$ for every $y\in Y$. For homogeneous polynomial, this
is equivalent to say that $A(y_1,\ldots,y_n)=0$ for all
$y_1,\ldots,y_n\in Y$, where $A$ is the symmetric multilinear form
associated to $P$. The study of zero subspaces of polynomials
depends heavily on whether the field of scalars is the field of real
numbers or the field complex numbers\footnote{For more information
about zero subspaces of polynomials on infinite dimensional Banach
spaces over $\mathbb{R},$ we refer to \cite{AroBoy} and the recent
work of Aron and H\'{a}jek~\cite{AroHaj}.}. For example, Aron and Hajek
\cite{AroHaj} has shown that every Banach space $X$ with weak-star
separable dual admits a homogeneous polynomial $P:X\To \mathbb{R}$
of degree 3 with no infinite-dimensional null space. It follows that
for example the Banach spaces $\ell_{\infty}$ and $\ell_1(\Gamma)$
for $\Gamma$ an index set of cardinality continuum admit homogeneous
polynomials of degree 3 with no infinite-dimensional null spaces. On
the other hand, for complex Banach spaces (on which we will
concentrate in this paper), Plichko and
Zagorodnyuk~\cite{PliZag} have proved the following result.\\

\begin{thm}\label{separablezero}
Let $X$ be an infinite-dimensional complex Banach space and let
$P:X\To\mathbb{C}$ be a polynomial with $P(0)=0$. Then, there is
an infinite-dimensional zero subspace $Y\subset X$ for $P$.\\
\end{thm}

\noindent As indicated above, in this paper we are mainly concerned
with polynomials with values in the field of complex numbers, so
from now on all Banach spaces that we will be working with will be
assumed to be over the field of complex numbers. Note that
Theorem~\ref{separablezero} suggests several natural questions
concerning zero subspaces of polynomials when the Banach space $X$
is assumed to be nonseparable, and not just infinite-dimensional.
Most prominent of these are the following two problems.\footnote{
Reiterated, for example, in the recent lecture by R. Aron at the
Caceres Conference on Banach Space Theory of
September 2006.}\\

\begin{problem}\label{first}
Let $X$ be a nonseparable Banach space and $P:X\To \mathbb{C}$ a
polynomial with $P(0)=0$. Is there always a nonseparable zero
subspace for $P$?
\end{problem}

\begin{problem}\label{second}
Let $X$ be a nonseparable Banach space and $P:X\To \mathbb{C}$ a
polynomial with $P(0)=0$. May $P$ have a separable maximal zero
subspace? Do all maximal zero
subspaces of $P$ have the same density?\\
\end{problem}

In this note, we shall show that the answer to both questions is in
general negative by providing the respective counterexamples as
quadratic functionals on the complex version of the Banach space
$\ell_1(\omega_1)$. It should be noted that the first problem has
already some history behind. Banakh, Plichko and
Zagorodnyuk~\cite{BanPliZag} noticed that if $X^\ast$ is weak$^\ast$
nonseparable, then every quadratic functional on $X$ has a
nonseparable zero subspace, and this is the case for instance of
$\ell_1(\mathfrak c^+)$ and of any nonseparable reflexive space.
Using this fact, Fern\'{a}ndez-Unzueta \cite{FerUnz} proved that
also every quadratic functional on $\ell_\infty$ has a nonseparable
zero subspace. Generalizing that result of~\cite{BanPliZag}, we
notice that if the weak$^\ast$-density of $X^\ast$ is greater than
$\beth_{m-2}$, then every homogeneous polynomial $P:X\To\mathbb{C}$
of degree $m$ has a nonseparable zero subspace. In this regard, we
may pose the
following problem:\\

\begin{problem}Let $m>1$. Is there a homogeneous polynomial
$P:\ell_1(\beth_{m-1})\To \mathbb{C}$ of degree $m$ such that all
its zero subspaces are
separable?\\
\end{problem}

 They explicitly
ask in~\cite{BanPliZag} whether there is a quadratic functional
$P:\ell_1(\omega_1)\To\mathbb{C}$ with all zero subspaces
separable and this is precisely what we construct. The
construction of the counterexample to problem \ref{second} is
rather standard and uses large almost disjoint families of sets.
On the other hand, the counterexample for question \ref{first}
requires the existence certain partitions $f:[\Gamma]^2\rightarrow
\{0,1\}$ with special properties\footnote{It should be noted
however that use of partitions of this kind for constructing
polynomials appears first in the paper of Hajek and Aron
\cite{AroHaj}.}. The proof of the existence of such partitions
relies on the machinery of minimal walks on countable ordinals
from~\cite{partitioning}. These partitions can also be produced in
a forcing extensions of the universe of sets, showing, for
example, the consistency of the existence of a quadratic
functional with no nonseparable zero subspaces exists over spaces
of the form $X=\ell_1(\Gamma)$ for $\Gamma$ an index set of
cardinality continuum. This leads us to the following natural
question to be answered if possible without appealing to special
axioms of set theory.

\begin{problem}
Is there a polynomial $P:\ell_1(\mathfrak c)\To \mathbb{C}$ with
$P(0)=0$ with no nonseparable zero subspace?
\end{problem}

Using another related special property of partitions arising again
from the work of the second author ~\cite{partitioning}, we produce
some vector-valued polynomials with even stronger properties. For
example, we show that for every separable Banach space $V$ there
exists a quadratic functional $P:\ell_1(\omega_1)\To V$ such that
the closure of the image under $P$ of any nonseparable subspace of
$\ell_1(\omega_1)$ has nonempty interior. Asking for this kind of
extension of our basic example is motivated by the work of
Baumgartener and Spinas~\cite{BauSpi} which also uses the partitions
of ~\cite{partitioning} to produce complex bilinear forms in the
context of vector spaces over a countable field rather than in the
context of Banach spaces over the field $\mathbb{C}$

\section{General remarks about polynomials}

The examples that we construct are defined on spaces
$\ell_1(\Gamma)$. Over this space, all polynomials can be
explicitly described. For instance, the general form of a
quadratic functional $P:\ell_1(\Gamma)\To\mathbb{C}$ is:

$$P(x) = \sum_{\alpha,\beta\in\Gamma}\lambda_{\alpha\beta}x_\alpha
x_\beta,\ \ \ \ x=(x_\gamma)_{\gamma\in\Gamma}\in\ell_1(\Gamma),$$

where $(\lambda_{\alpha\beta})_{\alpha,\beta\in\Gamma}$ is a
bounded family of complex scalars. Indeed in all our examples the
coefficients $\lambda_{\alpha\beta}$ are either 0 or 1, they
functionals of the form:

$$P(x)= \sum_{\{\alpha,\beta\}\in G}x_\alpha x_\beta,$$

where $G$ a certain set of couples of elements of $\Gamma$. On the
other hand, we state the following elementary basic fact about
polynomials that we shall explicitly use at some point:\\

\begin{prop}\label{uniformpolynomial}
Let $P:X\To Y$ be a homogeneous polynomial of degree $n$ and norm
$K$.
Then $\|P(x)-P(y)\|\leq n K M^{n-1}\|x\|\|y\|$ for every $x,y\in X$.\\
\end{prop}

\section{Polynomials with both separable and nonseparable maximal
zero subspaces}

Let $\Omega$ be a set and $\mathcal{A}$ be an almost disjoint
family of subsets of $\Omega$ (that is, $|A\cap A'|<|\Omega|$
whenever $A,A'\in\mathcal{A}$ are different), and let $\mathcal{B}
= \Omega\cup\mathcal{A}$. We consider the following quadratic
functional $P:\ell_1(\mathcal{B})\To\mathbb{C}$ given by
$$P(x) = \sum\{x_n x_A: n\in\Omega, A\in\mathcal{A}, n\in A\}.$$

\begin{thm}\label{smallmaximal}
The space $X=\ell_1(\Omega)\subset \ell_1(\mathcal{B})$ is
maximal zero subspace for the polynomial $P$.\\
\end{thm}

\emph{Proof}: The only point which requires explanation is that
$X$ is indeed maximal. So assume by contradiction that there is a
vector $y$ out of $X$ such that $Y=span(X\cup\{y\})$ is a zero
subspace for $P$. Without loss of generality, we suppose that $y$
is supported in $\mathcal{A}$. Pick $A\in \mathcal{A}$ such that
$|y_A| = max\{|y_B| : B\in\mathcal{A}\}$ and
$\mathcal{F}\subset\mathcal{A}$ a finite subset of $\mathcal{A}$
such that
$\sum_{B\in\mathcal{A}\setminus\mathcal{F}}|y_B|<\frac{1}{9}|y_A|$.
Now, because $\mathcal{A}$ is an almost disjoint family of subsets
of $\Omega$, it is possible to find $n\in\Omega$ such that $n\in
A$ but $n\not\in B$ whenever $B\in\mathcal{F}\setminus\{A\}$.
Consider the element $y+1_n\in Y$. We claim that $P(y+1_n)\neq 0$
getting thus a contradiction.

$$P(y+1_n) = \sum_{n\in B}y_B = y_A +
\sum_{n\in B, B\in\mathcal{A}\setminus\mathcal{F}}y_B$$

The second term of the sum has modulus less than $\frac{1}{9}$ the
modulus of the first term. So $P(y+1_n)\neq 0$.$\qed$\\

We are interested in the case when $|\mathcal{A}|>|\Omega|$. The
subspace $\ell_1(\mathcal{A})$ is a zero subspace for $P$. It may
not be maximal but this does not matter because, by a Zorn's Lemma
argument, it is contained in some maximal zero subspace. This fact
together with Theorem~\ref{smallmaximal} shows that $P$ has
maximal zero subspaces of both densities $|\Omega|$ and
$|\mathcal{A}|$.\\

There are two standard constructions of big almost disjoin
families. One is by induction, and it shows that for every
cardinal $\kappa$ we can find an almost disjoint family of
cardinality $\kappa^+$ on a set of cardinality $\kappa$. The other
one is by considering the branches of the tree $\kappa^{<\omega}$,
and this indicates that for every cardinal $\kappa$ we can find an
almost disjoint family of cardinality $\kappa^\omega$
(One construction or the other provides a better result depending on
whether $\kappa^\omega = \kappa$, $\kappa^\omega = \kappa^+$ or $\kappa^\omega>\kappa^+$). Hence,\\

\begin{cor}
Let $\kappa$ be an infinite cardinal and
$\tau=\max(\kappa^+,\kappa^\omega)$. There exists a quadratic
functional on $\ell_1(\tau)$ with a maximal zero subspace of
density $\kappa$ and another maximal zero subspace of density
$\tau$.\\
\end{cor}

\begin{cor}
There exists a quadratic functional on $\ell_1(\mathfrak c)$ with
a separable
maximal zero subspace and a maximal zero subspace of density $\mathfrak c$.\\
\end{cor}

\section{Polynomials for which all zero subspaces are separable}

We denote by $[A]^2$ the set of all unordered pairs of elements of
$A$,
$$[A]^2 = \{t\subset A : |t|=2\}.$$ We consider an ordinal $\alpha$
to be equal to the set of all ordinals less than $\alpha$, so
$$\omega_1 = \{\alpha : \alpha<\omega_1\}$$ is the set of countable
ordinals, and also for a non negative integer $n\in\mathbb{N}$,
$$n=\{0,1,\ldots,n-1\}.$$ We introduce some notations for subsets of a
well ordered set $\Gamma$. If $a\subset \Gamma$ is a set of
cardinality $n$, and $k<n$ we denote by $a(k)$ the $k+1$-th element
of $a$ according to the well order of $\Gamma$, so that
$$a=\{a(0),\ldots,a(n-1)\}.$$ Moreover, for $a,b\subset\Gamma,$ we write
$a<b$ if $\alpha<\beta$ for every $\alpha\in a$ and every
$\beta\in b$.\\

We recall also that a $\Delta$-system with root $a$ is a family of
sets such that the intersection of every two different elements of
the family equals $a$. The well known $\Delta$-system lemma
asserts that every uncountable family of finite sets has an
uncountable subfamily which forms a $\Delta$-system (see, e.g., \cite{EHMR}).\\

\begin{defn}\label{suitableColoring}
A function $f:[\Gamma]^2\To 2$ is said to be a \emph{partition of
the first kind} if for every uncountable family $A$ of disjoint
subsets of $\Gamma$ of some fixed finite cardinality $n$, and for
every $k\in n$ there exist $a,b,a',b'\in A$ such that
$f(a(k),b(k))=1$, $f(a'(k),b'(k))=0$ and $f(a(i),b(j)) =
f(a'(i),b'(j))$ whenever $(i,j)\neq (k,k)$. Notice that, passing to
a further uncountable subfamily $A$, we can choose such $a<b$ such
that, in addition,
$f(a(i),a(j))=f(a'(i),a'(j))=f(b(i),b(j))=f(b'(i),b'(j))$ for all
$\{i,j\}\in [n]^2$.
%(and also, though this
%can be deduced, such that
%$f(a(i),a(j))=f(a'(i),a'(j))=f(b(i),b(j))=f(b'(i),b'(j))$ for all
%$\{i,j\}\in
%[n]^2$).\\
\end{defn}

A natural projection of the square-bracket operation of
\cite{partitioning} (see also \cite{StevoCoherent}) gives the
following.
\begin{thm}
For $\Gamma=\omega_1$ there is a partition $f:[\Gamma]^2\To 2$ of
the first kind.
\end{thm}

\begin{rem}
Note that if some index set $\Gamma$ admits a partition
$f:[\Gamma]^2\To 2$ of the first kind then its cardinality is not
more than the continuum. So one would naturally like to know if some
index set $\Gamma$ of cardinality continuum admits such a partition.
While at the moment we do not know the answer to this problem, we do
know that the natural forcing of finite approximations to
$f:[\Gamma]^2\To 2$ leads to a generic extension of the universe of
sets which admits a partition $f:[\Gamma]^2\To 2$ of the first kind
on some index set $\Gamma$ of cardinality continuum that is
considerably larger than $\aleph_1.$
\end{rem}

To a given a partition $f:[\Gamma]^2\To 2$, we associate the
polynomial on $\ell_1(\Gamma)$ given by
$$P_f(x) = \sum_{\{i,j\}\in f^{-1}(1)}x_i x_j.$$

\begin{thm}
If $f:[\Gamma]^2\To 2$ is a partition of the first king and if $Y$
is a subspace of $\ell_1(\Gamma)$ with $Y\subset
P^{-1}_f(0)$, then $Y$ is separable.\\
\end{thm}

\emph{Proof}: We simplify the notation by writing $P$ in place of
$P_f$ and we assume on the contrary that we have $Y\subset
P^{-1}(0)$ a nonseparable space. We can construct recursively a
transfinite sequence $(\beta^i : i<\omega_1)$ of different
elements of $\Gamma$ and vectors $(y^i : i<\omega_1)$ with $y^i\in
Y$, $\|y^i\|=1$ and $y^i_{\beta^i}\neq 0$. By passing to a cofinal
subsequence, we can assume that there exists $\varepsilon\in
(0,1)$ such that $|y^i_{\beta^i}|>\varepsilon$ for every
$i<\omega_1$. This $\varepsilon$ will be fixed from now on in the
proof. We can also suppose that the sequence $(\beta^i)$ is
increasing. Since a polynomial is uniformly continuous on any
bounded set (Proposition~\ref{uniformpolynomial}), we pick a
$\delta$ such that $|P(x)-P(y)|<\frac{\varepsilon^2}{8}$ whenever
$\|x-y\|<\delta$ and $\|x\|,\|y\|\leq 3$.\\

For every $i<\omega_1$ we fix a finite set $a_i\subset supp(y^i)$
such that $\beta^i\in a_i$ and $$\sum_{u\not\in
a_i}|y^i_u|<\delta/3.$$ In other words, if we call $x^i$ to the
vector which coincides with $y^i$ on the coordinates of $a_i$ and is
supported on $a_i$, we have that $\|x^i-y^i\|<\delta/3$. By virtue
of the $\Delta$-system lemma, passing to a further cofinal
subsequence we can suppose that the sets $\{a_i : i<\omega_1\}$ form
a $\Delta$-system with root $a$. We call $b_i = a_i\setminus a$.
Again, by passing to a cofinal subsequence, we may suppose that:\\

\begin{itemize}
\item There exists $n<\omega$ such that $|b_i| =n$ for every
$i<\omega_1$.\\

%\item $a<b_i<b_j$ whenever $i<j$.\\

\item There exists $k\in n$ such $\beta^i = b_i(k)$ for every $i<\omega_1$.\\

\item There exists a function $e:a\times n\To 2$ such that
$f\{\alpha,b_i(r)\} = e(\alpha,r)$ for every $\alpha\in
a$ and $r\in n$.\\

\item There exists a function $d:[n]^2\To 2$ such that
$f\{b_i(r),b_i(s)\} = d\{r,s\}$ for every $\{r,s\}\in [n]^2$ and $i<\omega_1$.\\
\end{itemize}

%Thanks to part (1) of Theorem~\ref{suitableColoring}, we can
%suppose moreover that\\
%
%\begin{itemize}
%\item There exists a function $h:n^2\setminus\Delta_n\To \{R,G\}$
%such that for every $i<j$ and every $(r,s)\in
%n^2\setminus\Delta_n$,
%$f\{b_i(r),b_j(s)\} = h(r,s)$.\\
%\end{itemize}

For every possible function $q:n\times n\To 2$ we consider the
continuous function $\phi_q:\mathbb{C}^{a\cup
n}\times\mathbb{C}^{a\cup n}\To \mathbb{C}$ given as follows:

\begin{eqnarray*} \phi_q(z,z') &=& \sum\{(z_\alpha+z'_\alpha)(z_\beta+z'_\beta) :
\{\alpha,\beta\}\in [a]^2, f\{\alpha,\beta\}=1\}+\\
&& \sum\{(z_\alpha+z'_\alpha)(z_r+z'_r) : \alpha\in a, r\in n,
e(\alpha,r)=1\}+\\
&& \sum\{z_r z_s + z'_r z'_s : \{r,s\}\in [n]^2, d\{r,s\}=1\}+\\
&& \sum\{z_r z'_s : (r,s)\in n^2, q(r,s) = 1\}.
\end{eqnarray*}

That is, $\phi_q(z,z')$ computes the evaluation of the polynomial
$P$ on a hypothetical vector of the form $x^i + x^{i'}$ provided
that the restriction of $x^i$ to $a_i$ is given by $z$, the
restriction of $x^{i'}$ to $a_{i'}$ is given by $z'$, and the
value of $f(b_i(r),b_{i'}(s))$ is given by $q$. Let $\mathcal{B}$
be a countable base for the topology of the space
$\mathbb{C}^{a\cup n}$. For every $i<\omega_1$ we can choose
$U_i\in\mathcal{B}$ to be a neighborhood of the vector
$$y^i_\ast =
(y^i_{a(0)},\ldots,y^i_{a(|a|-1)},y^i_{b_i(0)},\ldots,y^i_{b_i(n-1)})$$
such that $|\phi_q(z)-\phi_q(w)|<\frac{\varepsilon^2}{8}$ for
every $z,w\in U_i\times U_i$ and every possible $q$. Since
$\mathcal{B}$ is countable, we can
make a further assumption:\\

\begin{itemize}
\item There exists $U\in\mathcal{B}$ such that $U_i=U$ for every
$i<\omega_1$.\\
\end{itemize}

Now, we use the fundamental property of $f$ to find ordinals $i<j$
and $i'<j'$ such that $f\{b_i(k),b_j(k)\}=1$,
$f\{b_{i'}(k),b_{j'}(k)\}=0$ and in all other cases,
$f\{b_{i}(r),b_{j}(s)\}=f\{b_{i'}(r),b_{j'}(s)\}$.\\

We claim that $P(y^i+y^j)\neq P(y^{i'}+y^{j'})$, which finishes
the proof since these two vectors belong to $Y$, and this way they
cannot be both zero.\\

Notice that
$|P(x^{i'}+x^{j'})-P(y^{i'}+y^{j'})|<\frac{\varepsilon^2}{8}$ and
$|P(x^{i}+x^{j})-P(y^{i}+y^{j})|<\frac{\varepsilon^2}{8}$ since we
know that $\|x^u-y^u\|<\delta/3$ for every $u$ and also
$|P(x)-P(y)|<\frac{\varepsilon^2}{8}$ whenever
$\|x-y\|<\delta$ and $\|x\|,\|y\|\leq 3$.\\

Define a function $q:n^2\To 2$ as $q(r,s) =
f\{b_{i'}(r),b_{j'}(s)\}$. We can compute\begin{eqnarray*}
P(x^{i'}+x^{j'}) &=& \phi_q(y^{i'}_\ast,y^{j'}_\ast),\\
P(x^i+x^j) &=& \phi_q(y^i_\ast,y^j_\ast) + y^i_{b_i(k)}y^j_{b_j(k)}.\\
\end{eqnarray*}

On the one hand, $y^i_{b_i(k)}y^j_{b_j(k)} =
y^i_{\beta^i}y^j_{\beta^j}$, so
 $|y^i_{b_i(k)}y^j_{b_j(k)}|>\varepsilon^2$.\\

On the other hand, since
$y^{i'}_\ast,y^{j'}_\ast,y^i_\ast,y^j_\ast\in U$, we have that
$$|\phi_q(y^i_\ast,y^j_\ast)-\phi_q(y^{i'}_\ast,y^{j'}_\ast)|<\frac{\varepsilon^2}{8}.$$

All these inequalities together yield that
$|P(y^{i'}+y^{j'})-P(y^{i}+y^{j})|>\frac{\varepsilon^2}{2}$.$\qed$\\

\section{Polynomial functionals with ranges in separable Banach spaces}

We shall denote by $$\Delta_n = \{(i,i) : i\in n\}$$ the diagonal of
the cartesian product $n\times n$.\footnote{Recall our convention,
$n=\{0,1,...,n-1\}.$} Also $B_X(x,r)$ or simply $B(x,r)$ will denote
the ball of center $x$ and radius $r$ in a given Banach space $X$.\\

\begin{defn}\label{omegaColoring}
 A function $f:[\Gamma]^2\To\omega$ is said to be a \emph{partition of the second kind} if for every uncountable
family $A$ of finite subsets of $\Gamma$ all of some fixed
cardinality
$n$, we have the following two conclusions\\
\begin{itemize}
\item [(a)] there is an uncountable subfamily $B$ of $A$ and a function
$h:n^2\setminus\Delta_n\To \omega$ such that $f(a(i),b(j)) =
h(i,j)$ for every $i\neq j$, $i,j<n$ and every $a<b$ in $B$;\\

\item [(b)] for every function $h:n\To\omega$ there exists $a<b$
in $A$
such that $f(a(i),b(i))= h(i)$.\\
\end{itemize}
\end{defn}

A natural projection of the square-bracket operation of
\cite{partitioning} leads again to the following result (see also
\cite{StevoCoherent}).

\begin{thm}
For $\Gamma=\omega_1$ there is a partition $f:[\Gamma]^2\To 2$ of
the second kind as well.
\end{thm}

\begin{rem}
Clearly, every such a set must have a cardinality no bigger than the
continuum. Similarly to the case of partitions of the first kind, we
know that it is consistent to have partitions of the second king on
some index set $\Gamma$ of cardinality continuum that is bigger than
$\aleph_1$ though this cannot be achieved by going to a generic
extension of the poset of finite approximations to such a partition.
\end{rem}

Fix a separable space $V$ and $\{v_n : n<\omega\}$ a dense subset of
the unit ball of $V$ with $v_0=0$. Then to every partition
$f:[\Gamma]^2\To\omega$ we associate a homogeneous polynomial
$P_f:\ell_1(\Gamma)\To V$ of degree 2 by the formula

$$P_f(x) = \sum_{\{i,j\}\in [\Gamma]^2} x_i x_j v_{f(i,j)}$$

\begin{thm}
Suppose that $f:[\omega_1]^2\To \omega$ is a partition of the
second kind and let $P=P_f:\ell_1(\Gamma)\To V$ be the
corresponding polynomial. Let $Y$ be a nonseparable subspace of
$\ell_1(\omega_1)$. Then
$\overline{P(Y)}$ has nonempty interior in $V$.\\
\end{thm}
\emph{Proof}: Let us introduce some language. Given $G,H\subset V$
and $\lambda>0$, we say that $G$ is $\lambda$-dense in $H$ if for
every $h\in H$ there exists $g\in G$
such that $\|g-h\|<\lambda$.\\

We can construct recursively a transfinite sequence $(\beta^i :
i<\omega_1)$ of different elements of $\Gamma$ and vectors $(y^i :
i<\omega_1)$ with $y^i\in Y$, $\|y^i\|=1$ and $y^i_{\beta^i}\neq
0$. By passing to a cofinal subsequence, we can assume that there
exists $\varepsilon\in (0,1)$ such that
$|y^i_{\beta^i}|>\varepsilon$ for every $i<\omega_1$. This
$\varepsilon$ will be fixed from now on in the proof. We can also
suppose that the sequence $(\beta^i)$ is increasing. Set
$\varepsilon_m = 2^{-m}\varepsilon^2/8$. Recursively on $m$, we
shall define a sequence $S_1\supset S_2\supset\cdots
S_m\supset\cdots$ of uncountable subsets of $\omega_1$ and a
convergent sequence $(w_m)$ of elements of $V$ such that for every
$m$, $P(Y)$ is $\varepsilon_m$-dense in $B_V(w_m,\varepsilon^2)$.
After this, $P(Y)$ will be dense in $B_V(\lim
w_n,\varepsilon^2/2)$ and the proof will be concluded.\\

In each inductive step $m$ we shall produce the uncountable set
$S_m\subset\omega_1$, and also for every $i\in S_m$ a finite set
$a_i^{(m)}\subset {\rm supp}(y^i)$ satisfying certain properties. We
suppose that we carried out the construction up to the step $m-1$
and we exhibit how to make step $m$.\\

The product of complex numbers is uniformly continuous on bounded
sets, hence we can find a number $\eta_m>0$ such that whenever
$\xi,\xi',\zeta,\zeta'\in B_{\mathbb{C}}(0,1)$ are such that
$|\xi-\zeta|<\eta_m$ and $|\xi'-\zeta'|<\eta_m$, then
$|\xi\xi'-\zeta\zeta'|<\varepsilon_m/7$. Since
$\{\zeta\in\mathbb{C} : \varepsilon<|\zeta|\leq 1\}$ can be
covered by finitely many balls of radius $\eta_m$, we can find an
uncountable subsequence $S'\subset S_{m-1}$ and a complex number
$\zeta_m$, with $|\zeta_m|>\varepsilon$, such that
$y^i_{\beta_i}\in B(\zeta_m,\eta_m)$ for every $i\in S'$ (hence,
whenever $\xi,\xi'\in
B(\zeta_m,\eta_m)$ then $|\xi\xi'-\zeta_m^2|<\varepsilon_m/7$).\\

We take $\delta_m<\varepsilon_m/42$. We know by
Proposition~\ref{uniformpolynomial} that
$\|P(x)-P(y)\|<\varepsilon_m/7$ whenever $\|x-y\|<\delta_m$ and
$\|x\|,\|y\|\leq 3$, indeed for every homogeneous polynomial $P$ of degree 2 and norm 1.\\

For every $i$ we fix a finite set $a_i\subset {\rm supp}(y^i)$ such
that $\beta^i\in a_i$ and $$\sum_{u\not\in a_i}|y^i_u|<\delta_m/2.$$
In other words, if we call $x^i$ to the vector which coincides with
$y^i$ on the coordinates of $a_i$ and is supported on $a_i$, we have
that $\|x^i-y^i\|<\delta_m/2$. By virtue of the $\Delta_m$-system
lemma, passing to a further cofinal subsequence $S''\subset S'$ we
can suppose that the sets $\{a_i : i<\omega_1\}$ form a
$\Delta$-system with root $a$. We call $b_i = a_i\setminus a$. We
shall suppose that all $b_i$'s have the same cardinality $n$, and
also that $a<b_i<b_j$ whenever $i<j$. Since these finite sets are
those defined in step $m$ we may also call $a_i = a_i^{(m)}$,
$a=a^{(m)}$, $b_i=b_i^{(m)}$. If $m>1$, we shall require some
coherence properties of the $a^{(m)}_i$'s with respect to the
$a^{(m-1)}_i$'s. First, from the very beginning we take
$a_i^{(m)}\supset a_i^{(m-1)}$. Also, by passing to a further
cofinal subsequence, we suppose that the relative position of
$b_i^{(m-1)}$ inside $b_i^{(m)}$ is the same for every $i$, namely,
that if $b_i = \{b_i(0)\prec \ldots \prec b_i(n-1)\}$ then there
exists $X\subset n$ such that $b_i^{(m-1)} = \{b_i(r) : r\in
X\}$ for every $i\in S''$.\\

Again, by passing to a further uncountable subsequence, we can suppose that $S''\subset S_{m-1}$ has the following extra properties:\\

\begin{itemize}

\item There exists $k\in n$ such $\beta^i = b_i(k)$ for every $i\in S''$ (indeed, unless we are in the first step, this is already guaranteed).\\

\item There exists a function $e:a\times n\To 2$ such that
$f\{\alpha,b_i(r)\} = e(\alpha,r)$ for every $\alpha\in
a$, $r\in n$, $i\in S''$.\\

\item There exists a function $d:[n]^2\To 2$ such that
$f\{b_i(r),b_i(s)\} = d\{r,s\}$ for every $\{r,s\}\in [n]^2$ and $i\in S''$.\\

\item There exists a function $h:n\times n\setminus \Delta_n\To 2$
such that $f(b_i(r),b_j(s)) = h(r,s)$ whenever $i,j\in S_m$, $i<j$
and $r\neq s$. We can do this since $f$ is a partition of the
second kind.
\end{itemize}

We consider the continuous function $\phi:\mathbb{C}^{a\cup
n}\times\mathbb{C}^{a\cup n}\To V$ given as follows:

\begin{eqnarray*} \phi(z,z') &=& \sum\{(z_\alpha+z'_\alpha)(z_\beta+z'_\beta)v_{f(\alpha,\beta)} :
\{\alpha,\beta\}\in [a]^2\}+\\
&& \sum\{(z_\alpha+z'_\alpha)(z_r+z'_r)v_{e(\alpha,r)} : \alpha\in
a, r\in n\}+\\
&& \sum\{(z_r z_s + z'_r z'_s)v_{d(r,s)} : \{r,s\}\in [n]^2\}.\\
\end{eqnarray*}

That is, $\phi$ makes what would be the evaluation of the
polynomial $P$ on the sum of two vectors supported on two items of
the $\Delta$-system, except that it does not compute the summands
corresponding to couples $\{b_i(r),b_j(r)\}$ which are the only
ones which are possibly different for different couples $i<j$.\\

Let $\mathcal{B}$ be a countable base for the topology of the
space $\mathbb{C}^{a\cup n}$. For every $i$ we can choose
$U_i\in\mathcal{B}$ to be a neighborhood of the vector
$$y^i_\ast =
(y^i_{a(0)},\ldots,y^i_{a(|a|-1)},y^i_{b_i(0)},\ldots,y^i_{b_i(n-1)})$$
such that $|\phi(z)-\phi(w)|<\varepsilon_m/7$ for every $z,w\in
U_i\times U_i$. Since $\mathcal{B}$ is countable, we can
find an uncountable $S_m\subset S''$ such that:\\

\begin{itemize}
\item There exists $U\in\mathcal{B}$ such that $U_i=U$ for every
$i\in S_m$.\\
\end{itemize}

Now, we use that $f$ is a partition of the second kind to find
ordinals $\{i_p < j_p : p<\omega\}\subset S_m$ such that
$f\{b_{i_p}(k),b_{j_p}(k)\}=p$, and in all other cases, we have that
$f\{b_{i_p}(r),b_{j_p}(r)\}=0$.\\

By the definition of $\delta_m$, we have that
$\|P(x^{i_p}+x^{j_p})-P(y^{i_p}+y^{j_p})\|<\varepsilon_m/7$
since $\|x^i-y^i\|<\delta_m/2$ for every $i$.\\

We can compute\begin{eqnarray*}
P(x^{i_p}+x^{j_p}) &=& \phi(y^{i_p}_\ast,y^{j_p}_\ast) + y^{i_p}_{b_{i_p}(k)}y^{j_p}_{b_{j_p}(k)}v_p.\\
\end{eqnarray*}

On the one hand, $y^i_{b_i(k)}y^j_{b_j(k)} =
y^i_{\beta^i}y^j_{\beta^j}\in B(\zeta_m,\eta_m)\cdot
B(\zeta_m,\eta_m)$, so
$$|y^i_{b_i(k)}y^j_{b_j(k)}-\zeta_m^2|<\varepsilon_m/7,$$
by the definition of $\eta_m$.\\

On the other hand, since $y^{i_p}_\ast,y^{j_p}_\ast,\in U$, we
have that
$$|\phi(y^{i_p}_\ast,y^{j_p}_\ast)-\phi(y^{i_0}_\ast,y^{j_0}_\ast)|<\varepsilon_m/7.$$

Set $w_m=\phi(y^{i_0}_\ast,y^{j_0}_\ast)$. Using that the vectors
$\{v_p : p<\omega\}$ are dense in the ball of $V$, all this data
together indicate that the vectors $\{P(y^{i_p}+y^{j_p}) :
p<\omega\}$ are $(\varepsilon_m/7 + \varepsilon_m/7
+\varepsilon_m/7)$ - dense in $B_V(w_m,\zeta_m^2)\supset
B_V(w_m,\varepsilon^2)$.\\

 We still have to check one fact in
order to guarantee that the sequence of the $w_m$'s will be
convergent, namely that $\|w_m -
w_{m-1}\|<\varepsilon_{m-1}<2^{-m+1}$. Notice that $w_m =
P'(x^{i_0}+x^{j_0})$  where $P'$ is the polynomial obtained from
$P$ by deleting the summands corresponding to pairs of $\omega_1$
of the form $\{b_i(r),b_j(r)\}$ for every $i<j$ and every $r<n$.
Let us call $x^{i}_{(m-1)}$, $x^{j}_{(m-1)}$, $i_0(m-1)$,
$j_0(m-1)$, etc. the corresponding objects of the previous step.
Notice that $\|x^{i_0} - x^{i_0}_{(m-1)}\|<\delta_{m-1}/2$ and
$\|x^{j_0} - x^{j_0}_{(m-1)}\|<\delta_{m-1}/2$ so
$$\|P'(x^{i_0}+x^{j_0}) -
P'(x^{i_0}_{(m-1)}+x^{j_0}_{(m-1)})\|<\varepsilon_{m-1}/7.$$ On
the other hand, we have that the norm of the difference
$$\|P'(x^{i_0}_{(m-1)}+x^{j_0}_{(m-1)}) -
P'(x^{i_0(m-1)}_{(m-1)}+x^{j_0(m-1)}_{(m-1)})\|$$ is equal to
$$ \|\phi_{(m-1)}(y^{i_0}_\ast,y^{j_0}_\ast) -
\phi_{(m-1)}(y^{i_0(m-1)}_\ast,y^{j_0(m-1)}_\ast)\|<\varepsilon_{m-1}/7,$$
since all those vectors belong to $U_{(m-1)}$. Now, observe that
we also have the equality
$P'(x^{i_0(m-1)}_{(m-1)}+x^{j_0(m-1)}_{(m-1)})=w_{m-1}$ thanks to
the coherence properties of the sets $a_i^{(m)}$ with respect to
the $a_i^{(m-1)}$'s. This finishes the proof.$\qed$\\

\section{Results about the existence of large zero
subspaces}\label{existenceoflarge}

Let us recall that the sequence $\beth_n$ of cardinals is defined
recursively by,
$$\beth_0=\aleph_0\mbox{ and }\beth_{n+1}=2^{\beth_n}.$$

\begin{lem}
Let $n,m$ be natural numbers and $X$ a Banach space with the
property that $dens(X^\ast,w^\ast)>\beth_{n+m-1}$. Then every
family of $\leq \beth_n$ many homogeneous polynomials of degree
$\leq m$ have a common
zero subspace of density $\beth_{n}^+$ (indeed $\beth_{n+1}^+$ for $m\geq 2$).\\
\end{lem}

Proof: We proceed by induction on $m$, the case $m=1$ being clear.
Let $\{A_\xi : \xi<\beth_n\}$ be a family of symmetric $\leq
m$-linear mappings on $X$. We construct recursively a sequence
$\{x_\alpha : \alpha<\beth_{n+1}^+\}$ of vectors of $X$ as
follows: For every $\alpha<\beth_{n+1}^+$, the multilinear maps of
the form
$A_\xi(x_{\alpha_1},\ldots,x_{\alpha_k},\ast,\ldots,\ast)$, for
$\{\alpha_1,\ldots,\alpha_k\}$ a nonempty family of less than $m$
ordinals less than $\alpha$ and $\xi<\beth_n$, constitute a family
of $\leq \beth_{n+1}$ many $(\leq m-1)$-linear forms. By the
inductive hypothesis there is a common zero subspace for all of
them of density $\beth_{n+1}^+$, and in particular we can choose
$x_\alpha$ in such zero space with $x_\alpha\not\in\overline{\rm
span}\{x_\beta : \beta<\alpha\}$. The $x_\alpha$'s are almost null
for the $A_\xi$'s, the only evaluations which may not vanish are
the $A_\xi(x_\alpha,\ldots,x_\alpha)\in \mathbb{C}$. But since
$\alpha$ runs up to $\beth_{n+1}^+$ and $\xi$ runs only up to $\beth_n$, we can pass to a subsequence
of cardinality $\beth_{n+1}^+$of the $\alpha$'s where
$A_\xi(x_\alpha,\ldots,x_\alpha)= r_\xi$ does not depend on
$\alpha$. Let $\zeta$ be a complex number with $\zeta^m = -1$. The
span of the vectors $\{x_\alpha + \zeta x_{\alpha+1}\}$ for
$\alpha$ even is a common
zero subspace for all the multilinear forms $A_\xi$.$\qed$\\

\begin{cor}
Let $X$ be a Banach space such that
$dens(X^\ast,w^\ast)>\beth_{m-2}$ and $P$ a homogenous polynomial
of degree $m$ on $X$, then $P$ has a nonseparable zero subspace.\\
\end{cor}

Proof: Let $A$ be the symmetric multilinear form corresponding to
$P$. We construct inductively vectors $(x_\alpha :
\alpha<\omega_1)$ in the following way: For every
$\alpha<\omega_1$ we consider all polynomials of the form
$A(x_{\alpha_1},\ldots,x_{\alpha_k},\ast,\ldots,\ast)$ with
$\alpha_1,\ldots,\alpha_k<\alpha$. This is a family of $\beth_0$
many homogeneous polynomials of degree $\leq m-1$, so by the
previous lemma they have a common zero subspace of density
$\beth_1^+$. In particular we can choose
$x_\alpha\not\in\overline{\rm span}(x_\beta : \beta<\alpha)$ in
this common zero subspace. Again all evaluations of $A$ in the
vectors $x_\alpha$ are null except perhaps
$A(x_\alpha,\ldots,x_\alpha)$. After multiplication by suitable
scalars, we may suppose that this value
$A(x_\alpha,\ldots,x_\alpha)$ is constant independent of $\alpha$.
Let $\zeta$ be a complex number with $\zeta^m = -1$. The zero
subspace for $A$ is then the closed linear span of the
vectors $\{x_\alpha + \zeta x_{\alpha+1}\}$ for $\alpha$ even.$\qed$\\

\end{document}